\documentclass[runningheads,a4paper]{llncs}

\usepackage{amssymb}
\usepackage{graphicx}
\usepackage{subfigure}
\usepackage{epstopdf,enumerate}
\usepackage{latexsym,amsmath,amscd,amssymb}
\usepackage{url,framed,color}   
\newcommand{\keywords}[1]{\par\addvspace\baselineskip
\noindent\keywordname\enspace\ignorespaces#1}

\begin{document}

\mainmatter  

\title{Covariant un-reduction for curve matching}

\titlerunning{Covariant un-reduction for curve matching}

%
%
\author{A. Arnaudon\inst{1}, M. Castrill\'on L\'opez\inst{2}, D.D. Holm\inst{1}}%
%
\authorrunning{Arnaudon et al.}

\institute{Dept. of Mathematics\\
Imperial College London\\
SW7 2AZ London, UK
\and 
Dept. Geometr\'{\i}a y Topolog\'{\i}a, 
Facultad de Ciencias Matem\'{a}ticas \\
Universidad Complutense de Madrid\\
28040 Madrid, Spain}
%
%

\maketitle

\begin{abstract}
		The process of un-reduction, a sort of reversal of reduction by the Lie group symmetries of a variational problem, is explored in the setting of field theories. 
		This process is applied to the problem of curve matching in the plane, when the curves depend on more than one independent variable. This situation occurs in a variety of instances such as matching of surfaces or comparison of evolution between species. 
		A discussion of the appropriate Lagrangian involved in the variational principle is given, as well as some initial numerical investigations.
\keywords{shape analysis, curve matching, Lagrange-Poincar\'e reduction, covariant field theory}
\end{abstract}


\section{Introduction}
The idea of un-reduction was introduced in \cite{CoHo} for the purpose of using reparametrisation by the action of the group $G=\mathrm{Diff}(S^1)$ to improve resolution of selected features in dynamics and optimal control problems, particularly in matching closed planar curves, whose configuration space $Q$ is the space of embeddings $Q=\mathrm{Emb}(S^1,\mathbb R^2)$, by redistributing grid points in $S^1$ along the curve. The un-reduction process used in \cite{CoHo} was based on \emph{reconstruction}, the inverse of Lagrangian reduction by symmetry \cite{HoMaRa1998}, which relates the solutions on $TQ/G$ to solutions on $TQ$.  This version of the un-reduction process was applied in the outer metric setting in \cite{CoHo,Cotter2012}. In this setting, the deformation of the shape of a curve in $\mathrm{Emb}(S^1,\mathbb R^2)$ was applied to the \emph{embedding} space, $\mathbb R^2$, completely independently of any reparametrisation of the embedded space, $S^1$. 

In contrast, the un-reduction technique introduced in \cite{BEGBH} seeks a family of equations, called the un-reduced equations, on $TQ$, whose solution projects onto those of a set of Euler-Lagrange equations on $T(Q/G)$. Thus, the un-reduction process used in \cite{BEGBH} is distinct from the reconstruction process used in \cite{CoHo,Cotter2012}. Moreover, the un-reduction approach introduced in \cite{BEGBH} raised an important issue for numerical applications in curve matching by optimal control, since it intertwined the reduction and reconstruction processes. Namely, the geodesic distance between two curves should be independent of their parametrisations. In particular, the measure of the deformation of the shape of a curve should be independent of its parametrisation.

Here we address this issue for curve matching from the viewpoint of a reformulation of the un-reduction scheme of \cite{BEGBH}, which gives a framework in the setting of inner metric discussed in \cite{BBM}, in which the shape deformation is applied to the \emph{embedded} space, rather than the \emph{embedding} space. 

In the process of addressing this issue, we will also generalise the un-reduction scheme by formulating it as a covariant space-time field theory. 
This generalisation gives us the freedom to introduce additional independent `time' or `space' variables for the purpose of coordinating comparisons among shapes.
Introducing additional independent variables allows more flexibility in making shape comparisons than, for instance, the time warp approach of \cite{Dur2013}, which does not increase the number of independent variables. For example, one could imagine making comparisons of cylindrical surfaces by assembling closed curves resolved on two-dimensional slices. 
In this case, the additional space variable would be transverse to the slices, and one would make comparisons of surfaces as single entities, rather than comparing the evolution in time alone of each slice independently. 
In addition, the covariant field theory generalisation of the un-reduction framework in the inner metric setting could lead to a variety of other applications, a few of which are mentioned in Section \ref{app-sec}.

This work may be summarized as an extension of \cite{BEGBH} in the following three directions: 
\begin{enumerate}[(1)]
		\item We promote the un-reduction formulation of \cite{BEGBH} in classical mechanics to a covariant field theory by following the same reasoning. Namely, we derive the Lagrange-Poincar\'e reduction of the shape space (Section \ref{LP-reduc}) and un-reduce it by including an independent parametrisation (Section \ref{un-reduction}). 
		\item Instead of the curvature weighted metric used in \cite{BEGBH}, we implement Sobolev metrics, which avoid the issue of arbitrary small geodesic distances (Section \ref{UR-curve}).
	\item We finish by illustrating this approach and assessing its validity with a few numerical experiments in the classical mechanical setting (Section \ref{numerics}). 
\end{enumerate}
A complete exposition of covariant un-reduction containing the proofs and applications in other areas, such as in theoretical physics, can be found in \cite{ArCaHo}.  

The main topic of the present paper is the first point listed above. We shall focus our discussion on the covariant Lagrange-Poincar\'e (LP) reduction by symmetry in the context of curve matching.
The symmetry will be $G=\mathrm{Diff}(S^1)$, the diffeomorphism group which acts on the configuration space of planar curves, $Q=\mathrm{Emb}(S^1,\mathbb R^2)$, as a reparametrisation of $S^1$. 
LP reduction is a general method to deal with noncanonical reductions, in which the configuration space is not the symmetry group. LP reduction allows the explicit derivation of the dynamical equations on the quotient space, $Q/G$, which, in our case, is the space of shapes, $\mathrm{Emb}(S^1,\mathbb R^2)/\mathrm{Diff}(S^1)$. 
The field theoretic version of LP reduction that we will use here was developed in \cite{CaRa,GBRa}, based on the classical reduction theory introduced in \cite{CeMaRa2001}. 
We will use a simplified version of this theory applied directly to curve matching. 

The second improvement relative to \cite{BEGBH} on the list above involves the Riemannian metric used to derive the un-reduction equations. 
As recently pointed out in \cite{BBM,BBMM}, use of Sobolev inner metrics avoids the problem of having arbitrarily small geodesic distances between two curves. 
In addition, the un-reduction equations take a simpler form with the Sobolev metric than with the curvature weighted metrics used in \cite{BEGBH}. 

Finally, we tested the un-reduction approach in a few numerical experiments where we considered an initial value problem with reparametrisation. 
Even though we simply chose the forward Euler method for this initial value problem, without any further modifications, it still converged to the expected solution. The success of these simple numerical simulations motivates us to go further in future work to consider boundary value problems in the full field theoretical framework as explained in Section \ref{cov-matching}. 

\section{The geometry of curve matching}
We start by recalling some basic facts about the geometry of curve matching that we will use throughout the text. We refer to earlier works such as \cite{MM,BBM,BBMM,BEGBH} for more details.  

\subsection{Reduction structure.}
Let $Q= \mathrm{Emb}^+(S^1,\mathbb R^2)$ be the manifold of
positive oriented embeddings from $S^1$ to $\mathbb R^2$. Elements
in $Q$ are maps $c(\theta)\in \mathbb R^2$ for $\theta\in S^1$ and
elements in the tangent space $T_cQ$ are pairs $(c,v)$ with $c\in
\mathrm{Emb}^+(S^1,\mathbb R^2)$ and $u \in C^\infty ( S^1 ,
\mathbb R^2)$ a parametrized vector field along the curve $c$, thus
$TQ = Q \times C^\infty (S^1 , \mathbb{R}^2)$.

We then consider the group $G= \mathrm{Diff}^+(S^1)$ of orientation
preserving diffeomorphisms of $S^1$ and its Lie algebra $\mathfrak
g= \mathfrak X(S^1)$. The group $G$ acts on the right in
$\mathrm{Emb}^+(S^1,\mathbb R^2)$ as reparametrisation of curves
$c$ and the reduced space is the space of shapes in $\mathbb R^2$
\begin{align}
    \Sigma := \frac{Q}{G} = \frac{\mathrm{Emb}^+(S^1,\mathbb
    R^2)}{\mathrm{Diff}^+(S^1)}.
\end{align}
The reduction of the phase space $TQ$ by the group of symmetry $G$ is a complicated space but can be decomposed via the introduction of a principal connection $\mathcal{A}: TQ\to \mathfrak g$ as follows
\begin{align}
		\frac{TQ}{G} \xrightarrow{\quad \mathcal A\quad } T\Sigma \oplus \tilde{ \mathfrak g}= T\left (\frac{Q}{G}\right )\oplus \frac{Q\times \mathfrak g}{G}.
		\label{LP-reduce}
\end{align}
The space $T\Sigma$ is then the tangent space to the space of shapes, and the adjoint bundle $\tilde {\mathfrak g}$ will encode the parametrisation velocity of the curve. 
This space seems rather abstract but it corresponds to having a Lie algebra attached to each point of the base space $\Sigma$, which means that for each shape we have the freedom to attach an arbitrary parametrisation velocity. 
The construction of the connection is straightforward for curve matching. Namely, given the velocity $u\in T_c Q$, we consider its tangent and normal decomposition,
\begin{align}
	u(\theta)=v(\theta)\mathbf{t} (\theta) + h(\theta) \mathbf{n}(\theta )
	\,,
\end{align}
where $(\mathbf{t} , \mathbf{n})$ is the orthonormal Frenet frame
along $c$ and $v(\theta), h(\theta)\in \mathfrak X(S^1)$ are scalar functions along the curve, parametrised by $\theta\in S^1$.       
We clearly have that $v(\theta){\bf t}(\theta)$ is a
vector tangent to the orbits of $G=\mathrm{Diff}^+(S^1)$.
This decomposition defines the principal connection $\mathcal A:TQ\to \mathfrak g$, which, when applied to $u$, gives $\mathcal A(u)= v$, the reparametrisation velocity.  
The horizontal part of $u$ is then $h(\theta) {\bf n}(\theta )$
and we have a decomposition $TQ=HQ\oplus VQ$. 
We will also need the curvature of $\mathcal A$, defined as $\mathcal B :=d^\mathcal{A}\mathcal A = d\mathcal A + [\mathcal A,\mathcal A]$, but its exact form will not be needed here. 
In fact, we shall skip any technical details which are not directly useful for the present work, and refer the interested reader to \cite{CeMaRa2001,CaRa,GBRa} for the full discussion of this construction.

\subsection{Field theoretical structure.}
We can now extend this reduction structure by promoting the classical system to a covariant field theory, see \cite{CaPeRa,GBRa}.
In order to do this, we consider an open domain $N\subset\mathbb R\times \mathbb
R$ endowed with the Euclidean metric, the associated coordinates $(t,x)$ and a volume form $\mathbf v=dtdx$. 
For simplicity we only consider a two dimensional  space-time manifold $N$, but more dimensions can be added in a straightforward way. 
The tangent space $TQ$ is then promoted to the jet bundle in order to capture the space-time direction used to compute a tangent vector. 
In this simple setting, the jet bundle has a simple geometric meaning, given by $J^1(N,Q)\simeq T^*N\otimes TQ$ and a generic element will be written
\begin{align}
    j^1c = c_t(\theta)(t,x) dt + c_x(\theta)(t,x)dx,
\end{align}
that is, $c_t,c_x\in TQ$ are the derivatives of a map $c:N\to Q$ along $t$ and $x$ respectively.
The generalisation of the time derivative is the divergence operator, given in this simple case by $\mathrm{div}\,j^1c=\frac{d}{dt}c_t + \frac{d}{dx}c_x$. 

\subsection{Riemannian metrics} 
The choice of a convenient Riemannian metric on $Q$ which is invariant with respect to the action of $G$ is an interesting topic. See, for example, \cite{BBMM} and \cite{BBM} and the references therein for more discussion.
In our case, invariance under reparametrisation is achieved by considering arclength integrations and derivatives $ds=|c_\theta|d\theta$ and $D_\theta = \frac{1}{|c_\theta|}\partial_\theta$. 
However, the main difficulty lies in the geometrical properties of the metric, as we deal with an infinite dimensional space. 
The natural $L^2$ metric
\begin{align}
		g_{L^2}(u,v)=\int _{S^1} \langle u(\theta) , v(\theta )\rangle ds,
\end{align}
with $u,v \in T_c Q$ such that $\langle \cdot, \cdot \rangle$ is a dot product in $\mathbb R^2$, is not very useful as it can lead to arbitrarily small geodesic
distances in both $Q$ and $Q/G$, see \cite{MM}. 
The problem can be overcome in the shape space $Q/G$ by the metric
\begin{align}
    g_\kappa(u,v)=\int _{S^1} (1+A\kappa (\theta)^2) \langle u(\theta) ,v(\theta )\rangle ds,
    \label{weighted-metric}
\end{align}
where $A>0$ and $\kappa$ denotes the Frenet curvature of the curve $c$, defined as
\begin{align}
\kappa := (D_\theta D_\theta c)\cdot J(D_\theta c)
= D_\theta\mathbf{t}\cdot\mathbf{n}\,.
    \label{kappa-def}
\end{align}
In fact, the weighted metric in \eqref{weighted-metric} can still have arbitrarily small geodesic distance in $Q$ along the fibres of the fibration $Q\to Q/G$. A metric with a well
defined Riemannian distance in both $Q$ and $Q/G$ may be obtained by
adding higher order derivatives of $u$ and $v$ in a Sobolev-type
expression as
\begin{align}
		g_{H^1}(u,v)=\int _{S^1} \left(\langle u(\theta) , v(\theta ) \rangle + A^2\langle D_\theta u(\theta) , D_\theta v(\theta )\rangle \right) ds\, .
\label{sobolev-metric}
\end{align}
We can collect these three cases (as well as others, see \cite{BBMM}) by using 
\begin{align}
    g_{\mathcal P}(u,v) = \int_{S^1} \langle u(\theta),\mathcal P v(\theta)
    \rangle ,
    \label{metric-general}
\end{align}
for a convenient choice of a $G$-invariant self-adjoint
pseudo-differential operator $\mathcal P$ which can depend on the
curve and its derivatives. In particular, the operator for
\eqref{weighted-metric} is $\mathcal{P}=1+A\kappa ^2$ and for
\eqref{sobolev-metric} we have $\mathcal P = 1-A^2D_\theta^2$. 
One additional advantage of the operator associated to
\eqref{sobolev-metric} is that it does not depend on the curve,
whereas the operator for \eqref{weighted-metric} depends on the
curvature of the curve where it is evaluated. This represents a
great simplification in the expression of the un-reduced
equations. 

\section{Reduction and un-reduction}
We are now ready to perform reduction by symmetry from the space of embeddings to the shape space using the covariant Lagrange-Poincar\'e reduction. 

\subsection{Lagrange-Poincar\'e reduction}\label{LP-reduc}
Let's begin by recalling the original problem of curve matching. 
The matching problem is a boundary value problem in $\Sigma=Q/G$ with Lagrangian $l:T\Sigma\to \mathbb R$.
Hamilton's principle states that the Euler-Lagrange equations associated with this Lagrangian yield the solution which minimises the action functional given by $l$. 
In practice, the matching problem is solved by using a shooting method for determining the initial momentum such that the curve at the final time matches the target curve within some specified tolerance. 
Instead, we will start in the larger space $Q$, where the numerical experiments can be done easily and, as a first step, reduce this system such that we recover the Euler-Lagrange equations on $\Sigma$. 

We project the variational principle defined for $L$ from $J^1(N,Q)$ to its quotient
$J^1(N,Q)/G = J^1(N,\Sigma)\, \oplus\, ( T^*N \otimes \tilde{\mathfrak g})$, where $\tilde{\mathfrak g}=(Q\times \mathfrak g)/G$, as in equation \eqref{LP-reduce}. 
Critical solutions are maps $\sigma : N \to T^*N\otimes
\mathfrak{\tilde{g}}$ which, moreover, project to maps $\rho : N
\to \Sigma = Q/G$ as $\rho = \pi_{\mathfrak{\tilde{g}}} \circ
\sigma$ according to the diagram
\begin{equation}
\label{maps}
\begin{array}{lll}
&  & \hspace{-0.21in}T^{\ast }N\otimes \mathfrak{\tilde{g}} \\
\multicolumn{1}{c}{} &
\multicolumn{1}{c}{\overset{\hspace{-0.07in}\sigma } {
\nearrow }} & \downarrow ^{\pi _{\mathfrak{\tilde{g}}}} \\
\multicolumn{1}{c}{N} & \multicolumn{1}{c}{\overset{\rho
}{\longrightarrow }} & \Sigma
\end{array}
\end{equation}
where $\pi_{\mathfrak{\tilde{g}}}:
T^{\ast}N\otimes\mathfrak{\tilde{g}}\rightarrow \Sigma$ is the
projection of the adjoint bundle neglecting the $T^*N$ factor. The
free variations of the initial problem provide a family of
constrained variations that define a new type of variational
equations, called Lagrange-Poincar\'e equations, \cite{CaRa}, \cite{GBRa}. 
The next theorem gives the
Lagrange-Poincar\'e reduction which includes forces $F:T^*N\otimes TQ\to TQ$. 
Before stating this theorem without proof, we will make another important assumption which is satisfied by most of the Lagrangians used in curve matching. Namely, we assume our Lagrangian decomposes as a sum of two Lagrangians taking values from the vertical and horizontal space that will be denoted $L=L^h+L^v$ and $\ell=\ell^h+\ell^v$.

\begin{theorem}[Covariant Lagrange-Poincar\'{e} reduction with
forces]\label{LP}
    Given a map $c:N\rightarrow Q$, let $\sigma:N\rightarrow T^{\ast}N \otimes\mathfrak{\tilde{g}}$ be defined as 
\begin{align}
\sigma(x)(\omega )= [s(x),\mathcal A(Ts\cdot (\omega))]_G
\,,
\label{sigma-def}
\end{align}
	with $\omega\in T_xN, x\in N$ and where $[\cdot ]_G$ stands for the quotient by $G$; $\rho :N\to \Sigma$, $\rho(x) = [s(x)]_G = \pi_{\mathfrak{\tilde{g}}} \circ
    \sigma$. With the previous definitions, the following points are equivalent 
	\begin{enumerate}[(1)]
\item $s$ is a critical mapping of the variational principle
\begin{align}
    \delta\int_{N}L(s,j^{1}s)\mathbf{v}+\int_{N}\langle F(s,j^{1}s),\delta s\rangle\mathbf{v}=0
    \label{un-red-var}
\end{align}
with free variations $\delta s$.

\item The Euler-Lagrange form of $L$ satisfies the relation
\begin{align*}
    \mathcal{EL}(L\mathbf{v})(j^2 s)=F,
\end{align*}
where $\mathcal{EL}$ is the Euler-Lagrange operator acting on the second jet bundle (second order field theoretical tangent space), which gives the usual Euler-Lagrange equations. 
\item $\sigma:N\rightarrow T^{\ast}N\otimes$
$\mathfrak{\tilde{g}}$
is a critical mapping of the variational principle%
\begin{align*}
        \delta\int_{N}\ell(j^{1}\rho,\sigma)\mathbf{v}+\int_N  \langle f^{h}(j^{1}\rho,\sigma),
        \delta\rho\rangle\mathbf{v}+\int_N\langle f^{v}(j^{1}\rho,\sigma),\eta\rangle \mathbf{v}=0,
\end{align*}
for variations of the form
$\delta\sigma=\nabla\eta-[\sigma
,\eta]+\mathcal{B}(\delta\rho,T\rho)\in
\tilde{\mathfrak{g}}$, where $\delta\rho\in T_\rho \Sigma$ is a
free variation of $\rho$ and $\eta$ is a free section of
$\mathfrak{\tilde{g}}\rightarrow \Sigma$.

\item $\sigma$ satisfies the Lagrange-Poincar\'{e} equations, written if $L$ and $\ell$ decomposes according to the vertical/horizontal decomposition
\begin{equation}
\left.
\begin{array}{c} \mathrm{div}
\left(\dfrac{\delta \ell ^h}{\delta j^1 \rho}\right)
-\dfrac{\delta \ell ^h}{\delta \rho} = f^h + \dfrac{\delta \ell
^v}{\delta  \rho} - \left\langle \dfrac{\delta \ell ^v}{\delta
\sigma},\mathcal{B}(T\rho,\cdot)\right\rangle , \\
\mathrm{div}\dfrac{\delta\ell
^v}{\delta\sigma}+\mathrm{ad}_{\sigma}^{\ast}\dfrac{\delta\ell
^v}{\delta\sigma}=f^{v}.
\end{array}
\right\} \label{ff2}
\end{equation}
\end{enumerate}
\end{theorem}

One recognises left hand side of the first equation in \eqref{ff2} as an Euler-Lagrange equation and the second one as an Euler-Poincar\'e equation. 
The right hand side of both equations are either forces, or coupling with between them. 
The solution $\sigma$ of the Euler-Poincar\'e equation in \eqref{ff2} will influence the Euler-Lagrange equation via the term involving the curvature of the connection $\mathcal A$. An additional coupling arises because $\sigma $ is in the adjoint bundle and therefore depends implicitly on the base curve in $Q$.  

\subsection{Un-reduction}\label{un-reduction}
The particular form of the equations in \eqref{ff2}, based on the decomposition of the Lagrangian and the inclusion of the force term will allow us to decouple these equations in the sense that the right hand side of the EL equation will vanish; so the feedback of the EP equation to the EL equation will disappear. 
Before stating the un-reduction theorem we must recall that the canonical momentum map $\mathbf J:T^*Q \to \mathfrak{g}^*$ for the natural lift action of $G$ on $T^*Q$, is defined by 
\begin{align*}
\label{momapp}
\langle \mathbf J (\alpha_q),\xi\rangle_{\mathfrak g\times\mathfrak g^*}  = \langle \alpha_q,\xi_Q\rangle_{TQ\times TQ^*}\,,
\end{align*}
where $\alpha_q\in T^*Q$, $\xi\in \mathfrak g$, and $\xi_Q\in TQ$ is the infinitesimal transformation of the action of $G$ on $Q$ at the point $q\in Q$.
The map $\mathbf{J}$ extends to a map $\mathbf{J}: TN\otimes T^*Q \to TN\otimes \mathfrak{g}^*$, trivially in the factor $TN$.

We can finally state the covariant un-reduction theorem. We refer to \cite{ArCaHo} for the proof and more details about this theorem. 
\begin{theorem}
\label{mainth} 
We consider a $G$-equivariant force $F:J^1(N,Q)\to T^*Q$ such that
$F^v=p^v \circ F$ is arbitrary and $F^h = p^h \circ F$ is given by
the condition
\begin{equation}
\label{fh}
f^h = -\frac{\nabla \ell ^v}{\delta \rho} + \left\langle
\frac{\delta \ell ^v}{\delta
\sigma},\mathcal{B}(T\rho,\cdot)\right\rangle ,
\end{equation}
for its projection $f^h:J^1(N,\Sigma)\times
(T^*N\otimes\mathfrak{\tilde{g}})\to T^*\Sigma$. Then, the
variational equations of the problem defined by $L$ and $F$ read
\begin{equation}
\label{eqsp} \left.
\begin{array}{c} \mathcal{EL}(L^h)(j^2 s) =0\\
\mathcal{A}^* \mathrm{div}\left(\mathbf J\left(\dfrac{\delta
L^v}{\delta j^1 s}\right)\right)  =  F^v (j^1 s),
\end{array}
\right\}
\end{equation}
where $\mathcal{A}^* : \mathfrak{g}^*\to V^*Q$ is the dual of the
connection form. Finally, critical solutions $s:N\to Q$ of
\eqref{eqsp} project to critical solutions $\rho = [s]_G$ of the
Euler-Lagrange equations $\mathcal{EL}(l)(j^2 \rho)=0$.
\end{theorem}
\begin{remark}
    For $N=\mathbb R$, $\mathbf{v}=dt$, that is, in the case of classical mechanics,
    we have $\mathrm{div}=d/dt$ and we recover the
    results and equations of \cite{BEGBH}.
\end{remark}
The first equation in \eqref{eqsp} is the usual Euler-Lagrange equation for the horizontal Lagrangian $L^h$, needed for solving the matching problem.  
Regarding the interpretation of the second equation, the definition of $\mathbf{J}$ above shows
that $\mathbf{J}(\delta L^v/\delta j^1s)$ is a covariant momentum
map, so that $\mathrm{div}\, \mathbf{J}(\delta
L^v/\delta j^1s)$ is the expression of a conservation law with respect to the group of symmetries. 
If one set $F^v=0$, the conservation law is complete. However, sometimes it may be
interesting to keep this vertical force, as it might be used to
externally control the dynamics along the vertical space; that is, the reparametrisation.

\subsection{Un-reduction with Sobolev metric}\label{UR-curve}

We consider the $\mathrm{Diff}^+(S^1)$-invariant Lagrangian
$L:J^1(N,Q)\simeq T^*N\otimes TQ \to \mathbb{R}$ which can be decomposed as $L=L^h + L^v$ with respect to the connection $\mathcal{A}$ as
\begin{align*}
L^h(j^1c) = \frac12 \int_{S^1}\left ( h_t \mathcal Ph_t +  h_x \mathcal P h_x \right ) ds,\quad
L^v (j^1c) = \frac12 \int_{S^1}\left (v_t \mathcal Pv_t  + v_x \mathcal P v_x \right ) ds,
\end{align*}
where
\begin{equation*}
        c_t = v_t\mathbf{t} + h_t \mathbf{n} \quad \mathrm{and}\quad c_x= v_x\mathbf{t} + h_x\mathbf{n}.
\end{equation*}
The un-reduction equations \eqref{eqsp} are then computed in 
Proposition \ref{ppop} below in the case when $\mathcal P$ is
independent of the curve, that is, Sobolev metrics.

\begin{proposition}\label{ppop}
        The un-reduced equations \eqref{eqsp} for the bi-dimensional problem of planar
simple curves defined by the Lagrangian defined above and the metric \eqref{metric-general} with $\mathcal P$ being the Sobolev operator, read 
\begin{align}
        \begin{split}
				\partial _{x}\mathcal Ph_{x} + \partial _{t}\mathcal Ph_{t}&= D_\theta(h_{x}\mathcal P v_{x}+h_{t}\mathcal P v_{t})-\kappa H \\
    \partial_x \mathcal P v_x+\partial_t\mathcal P v_t&=F^v
        \end{split}
        \label{un-red-curve}
\end{align}
for any choice of vertical force $F^v$, where
\begin{align}
		H=\frac12\left (h_x \mathcal Ph_x+h_t\mathcal P h_t\right ).
\end{align}

\end{proposition}
In \eqref{un-red-curve}, the function $H$ is the shape kinetic energy density and $\kappa$ is the Frenet curvature of the curve $c$, defined in \eqref{kappa-def}.
This term can be interpreted as a penalty term in deforming curved regions. 
The sign of this term would depend on the concavity or convexity of the curve at this point, and thus this force would try to prevent the curve to be deformed too fast in these regions. 
Equation \eqref{un-red-curve} shows that the dynamics in $(x,t)$ is governed by the coupling between $h_t$ and $v_t$ required for the shape deformation to be independent of the reparametrisation.
In fact, it also contains a derivative with respect to $\theta$ and vertical vectors.
This is the only term which couples with the vertical equation, and it also gives the corrections necessary for the curve deformation to be independent of the reparametrisation.

\begin{remark}
		In the classical mechanics setting, the un-reduction equations with Sobolev metric would be very similar to \eqref{un-red-curve}, but without the $x$-dependent terms. 
		They will be the equations used for the numerics in section \ref{numerics}. 
		Owing to the simplicity of the Sobolev metric compared to the curvature weighted metric, the derivation of this equation is directly done from the un-reduction equations, not from the variational principle, as in \cite{BEGBH}. 
		We refer to \cite{ArCaHo} for the details of this calculation.  
\end{remark}

\section{Applications}\label{app-sec}

Before discussing the possible applications of the covariant un-reduction scheme in curve matching, we shall present a short numerical study of the un-reduction in classical mechanics using the Sobolev metric $H^1$. 

\subsection{Numerical validation}\label{numerics}

In order to test and illustrate the un-reduction scheme, we performed a few simple numerical experiments.
We restricted ourselves to the classical case, already done theoretically in \cite{BEGBH}, but with the $H^1$ norm instead of the curvature weighted norm. 
The only effect of the Sobolev norm that which interested us is that it regularises the curve deformation and prevents large bending of the curve, smaller that the scale given by $A$. 
A deeper analysis of the effect of the Sobolev norm in the matching process is not the aim of this paper but is important for applications.
Our main goal here was to check the decoupling between the shape and the reparametrisation dynamics for a simple initial value problem. 

\begin{figure}[htpb]
		\centering
		\subfigure[$dt=0.04$]{\includegraphics[scale=0.47]{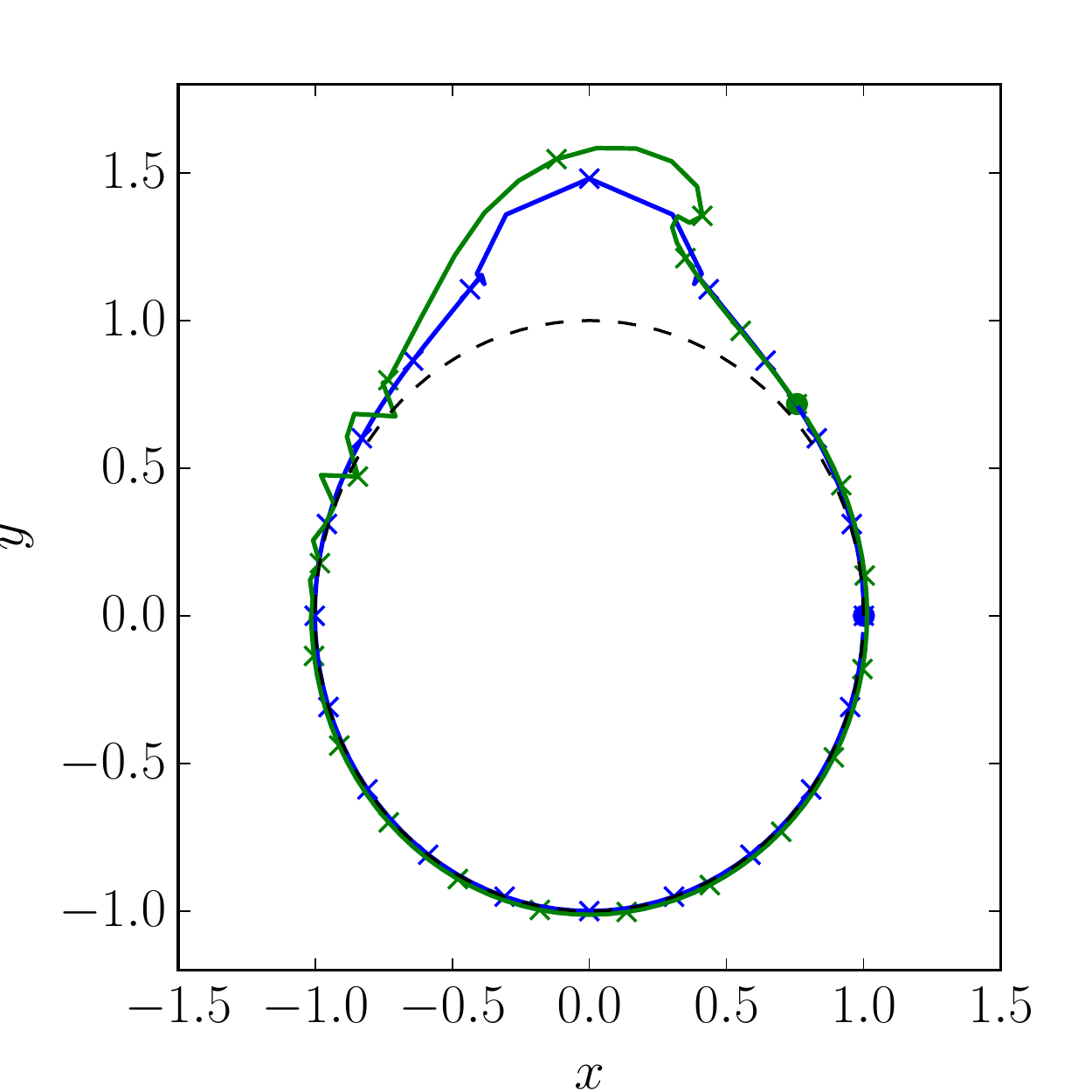}}
		\subfigure[$dt=0.02$]{\includegraphics[scale=0.47]{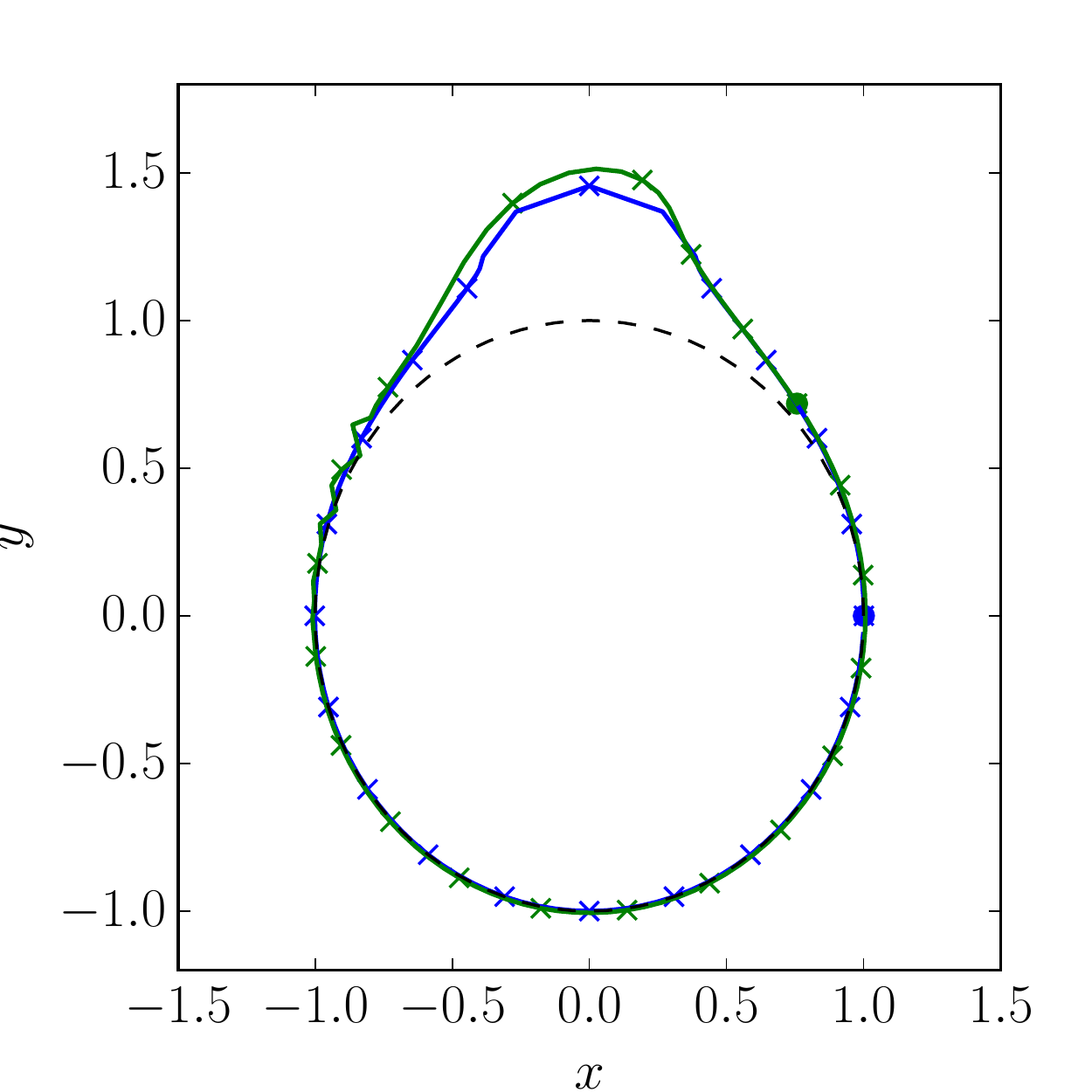}}
		\subfigure[$dt=0.01$]{\includegraphics[scale=0.47]{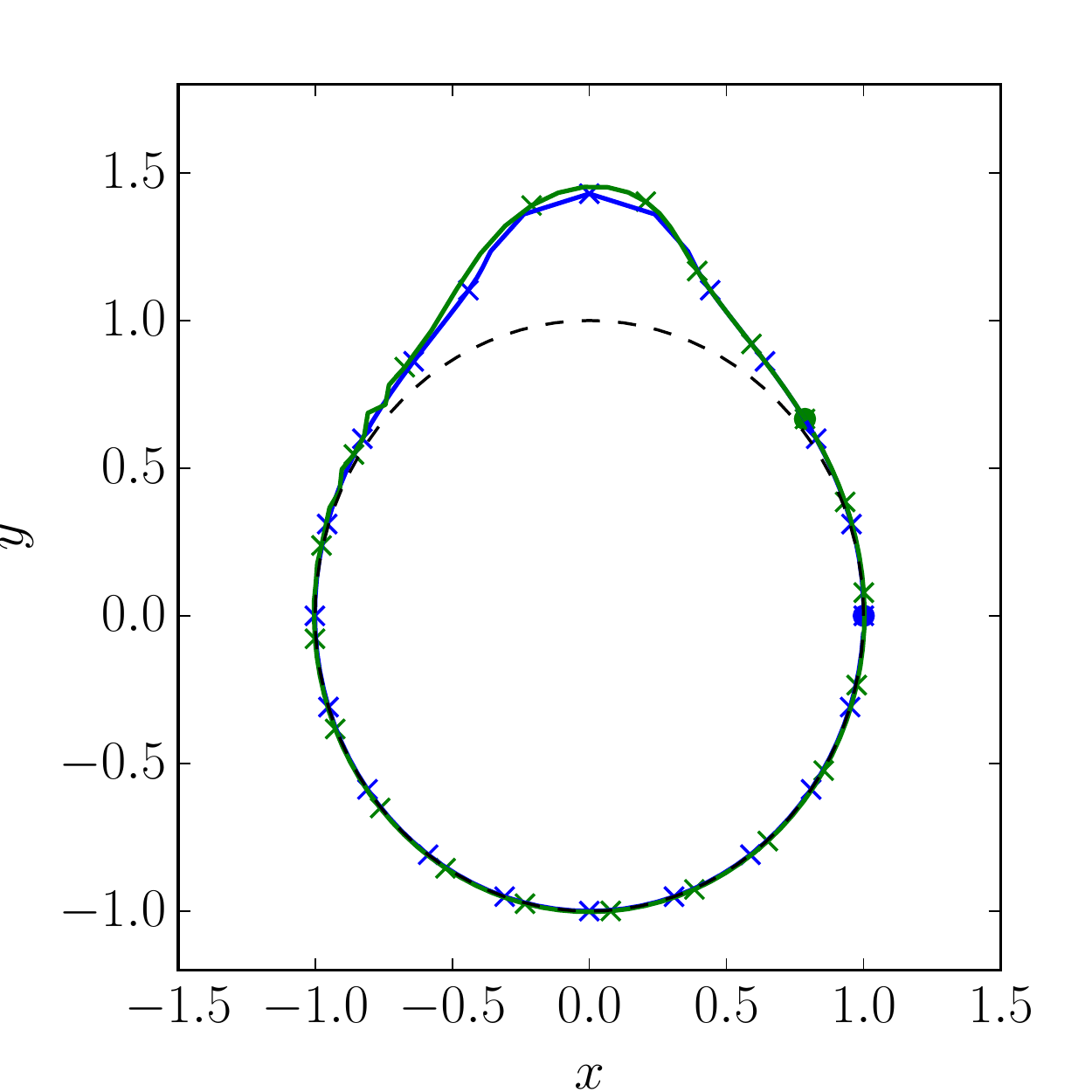}}
		\subfigure[Distance between fixed and reparametrised curves with the method of currents]{\includegraphics[scale=0.47]{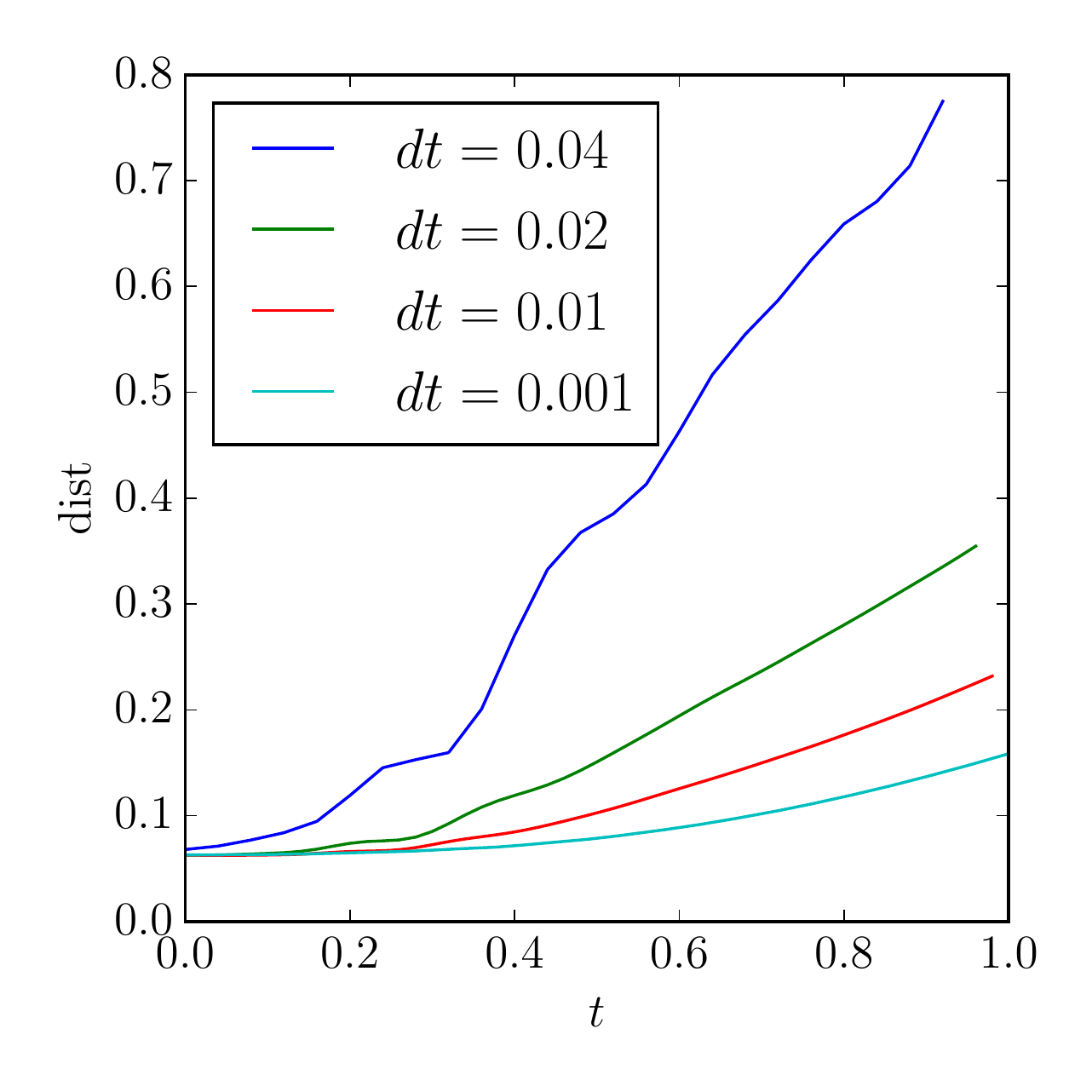}}
\caption{This figure illustrates the effect of the time step of $S^1$-reparameterisation on the quality of the deformation of a curve using the explicit Euler scheme in the un-reduction approach. 
The blue curve has horizontal initial velocity $h_t$ for shape deformation, while the green curve has an additional constant vertical initial velocity $v_t$, for reparameterisation. 
The blue and green dots begin at the same initial point, but then the green one shifts along the green curve as the reparameterisation proceeds. 
The black dashed curve is the initial condition of the simulation. 
Upon decreasing the time steps, the coupling between the vertical and horizontal dynamics decreases and the quality of the deformation improves, even when using an explicit Euler scheme. 
On the bottom right panel, we plotted the distance between the two curves, as a function different time resolutions. The distance is computed by the method of currents, which is independent of the parametrisation. 
}
		\label{fig:un-red}
\end{figure}

Our numerical scheme used the Euler explicit scheme in time and 2nd order centred finite difference approximation for $D_\theta$ in order to have a symmetric space discretisation in $\theta \leftrightarrow -\theta$. 
The application of the Sobolev operator $\mathcal P=1-A^2D_\theta^2$ was done in Fourier space with $A=0.3$.
Our initial condition was a circle and the initial horizontal velocity was a bump function, so the curve deformed as in Fig.~\ref{fig:un-red}. 
The curve was discretised with $100$ points and we used a set of time steps ($dt=0.04,0.02,0.01,0.001$) to study the convergence of the scheme, and especially the decoupling between the reparametrisation and the shape deformation.  
In order to do this, we ran two initial value problems where one of them also had a vertical initial velocity, taken to be constant such that the parametrisation would rotate during the evolution of the curve. 
We then computed the distance between the two curves using the methods of currents \cite{GTY2006} at each time to make a parametrisation-free comparison of the shapes of the curve. 
The results are displayed in Fig.~\ref{fig:un-red} together with the distance between the curves as a function of time. 
Even with the simple numerics we used (the Euler scheme and finite difference), the un-reduction feature was verified. 
The example we studied is simple and did not require high resolutions and robustness tests as for more realistic shapes. 
Further numerical studies using the un-reduction scheme would thus include improvements of the current implementation and a shooting method in order to solve the correct matching, or boundary value problem. 

\subsection{Covariant matching}\label{cov-matching}

In the simplest case of field theory treated here, namely a two dimensional space-time, two main applications for boundary value problems (BVP) present themselves for further discussion. 
Besides BVPs, initial value problems (IVP) could also be considered, but IVP are not of great interest for curve matching. We will thus forego discussing them here, although a possible application would be to predict the evolution of a particular model, knowing that it should ``roughly'' follow a generic model.  In this case, the initial value problem must have carefully chosen initial values, a subject which is out of the scope of this work. 
The two applications for BVP that we will discuss are the following:
\begin{enumerate}[(1)]
		\item Matching between cylindrical surfaces, and 
		\item Spatio-temporal analysis.
\end{enumerate}

{\bf Case (1)} In the first case, a set of slices along a cylindrical surface ( a typical example would be a bone) are given, where $x$ is the parameter along the main axis of the surface. 
For the sake of simplicity, we will just consider two slices, but more could be added without to much trouble. 
The first step would be to generate the initial and final conditions, namely use un-reduction for the IVP to interpolate between the two curves and generate the initial and final surfaces. 
Once this is done, covariant un-reduction could be applied using a shooting method in time, such that the solution is a critical point of the action functional $\int L(c,c_t,c_x) dtdx$. 
In our simple case, where the Lagrangian is purely quadratic, the solution would be a harmonic map, or a minimal surface, and would then require more advanced mathematics, beyond the present discussion. 
This model would compute the distance between two surfaces, taking into account that the interpolation between the slices in space should be imposed simultaneously with the matching in time. 
The resulting distance will be different than a naive model, which would compute the matching in time, slice by slice.  For an illustration of matching slice by slice, we refer to the last example in \cite{Cotter2012} where a surface representing a nasal cavity is reconstructed out of a set of slices. 
The step done there corresponds to the generation of initial and final surfaces, whereas covariant un-reduction would compute the distance between these two surfaces, during the temporal deformation of the entire nasal cavity. 

{\bf Case (2)} 
Spatio-temporal analysis has recently been reviewed in \cite{Dur2013}, from yet another viewpoint. 
Indeed, the matching in space done  in \cite{Dur2013} does not depend on a space parameter, but is \emph{instantaneous}, namely given by a single map between the two curves.
They also included a ``time warp'' which account for the change of pace of the evolution of the two models to be compared. 
In our case, the spatial variable comes into play on the same footing as time, and, thus, brings more flexibility into the comparison. 
Again, the theory of harmonic maps could help in understanding the properties of the solutions, and it is possible that the time warp reparameterisation could be recovered as well. \\

We finally want to mention the freedom to choose the vertical force in the un-reduction equations.
This force could be used to control the parametrisation during the matching. 
Different types of forces could be considered, such as a force which would optimally redistribute the parametrisation in different regions of the curve, such that the number of points for discretising the curve could be reduced.
Another force could be used to match the paramerisation of the target curve, such that the computationally more expensive method of currents could be avoided.  

\subsubsection*{Acknowledgements}
We are grateful to M. Bauer, M. Bruveris, L. Younes, S. Durrleman, R. Montgomery and T.
Ratiu for valuable discussions during the course of this work. AA
acknowledges partial support from an Imperial College London Roth
Award, AA and DH from the European Research Council Advanced Grant
267382 FCCA. MCL has been partially funded by MINECO (Spain) under
projects MTM2011-22528 and MTM2010-19111. MCL wants to thank
Imperial College for its hospitality during the visit in which the
main ideas of this work were developed.


\begin{thebibliography}{99}

		\bibitem{ArCaHo} A. Arnaudon, M. Castrill\'on L\'opez, D.D. Holm,
				\emph{Covariant un-reduction and applications}, in preparation.

		\bibitem{BBMM} M. Bauer, M. Bruveris, S. Marsland, P.W. Michor,
				\emph{Constructing reparameterization invariant metrics on spaces
				of plane curves}, Differential Geom. Appl. {\bf 34} (2014),
				139--165.

		\bibitem{BBM} M. Bauer, M. Bruveris, P.W. Michor, \emph{Overview of
				the Geometries of Shape Spaces and Diffeomorphism Groups}, J.
				Math. Imaging and Vision {\bf 50} (2014), 67--90.

		\bibitem{BEGBH} M. Bruveris, D.C.P. Ellis,  F. Gay-Balmaz, D.D. Holm,
				\emph{Un-reduction}, Journal of Geometric Mechanics {\bf 3} (2011), 363--387.

		\bibitem{CaRa} M. Castrill\'on L\'opez, T.S. Ratiu, \emph{Reduction in Principal Bundles:
				Covariant Lagrange-Poincar\'e Equations}, Comm. Math. Phys. {\bf 236} (2003), 223-–250.

		\bibitem{CaPeRa} M. Castrill{\'o}n L{\'o}pez, P.L. Garc\'\i a, T.S. Ratiu,
				\emph{ Euler--Poincar{\'e} reduction on principal bundles}, Lett. Math. Phys. {\bf 58} (2001), 167--180.

		\bibitem{CeMaRa2001} Cendra, Hern{\'a}n and Marsden, Jerrold E and Ratiu, Tudor S
				\emph{Lagrangian reduction by stages}, American Mathematical Soc.
				{\bf 722}, (2001), x+108 

		\bibitem{Cotter2012} C.J. Cotter, A. Clark, J. Peir{\'o}, \emph{A reparameterisation based approach to geodesic constrained solvers for curve matching}, Int. J. Comput. Vis.
				{\bf 99}, (2012), 103¿121

		\bibitem{CoHo}  C.J.  Cotter  and  D.D.  Holm, \emph{Geodesic  boundary  value
				problems  with  symmetry},  J.  Geom. Mech. {\bf 2}, no. 1 (2010),
				417--444.

		\bibitem{Dur2013} S. Durrleman, X. Pennec, A. Trouv\'e, J. Braga, G. Gerig, N. Ayache, 
				\emph{Toward a comprehensive framework for the spatiotemporal statistical
				analysis of longitudinal shape data}, Int J Comput Vis. {\bf 103},
				no. 1 (2013), 22--59.

		\bibitem{GBRa} D.C. Ellis, F. Gay-Balmaz, D.D. Holm, T.S. Ratiu,
				\emph{Lagrange-Poincar\'e field equations}, Journal of Geometry
				and Physics {\bf 61}, no. 11, (2011), 2120--2146.

		\bibitem{GTY2006} J. Glaunes, A. Trouve, L. Younes,
			\emph{Modeling planar shape variation via hamiltonian flows of curves}, in H. Krim Jr. and Y.A (Eds.), Statistics and analysis of shapes (2006), Basel Birkh\"auser.



		\bibitem{HoMaRa1998} 
				D.D. Holm, J.E. Marsden and T.S. Ratiu,
				The Euler--Poincar\'e equations and semidirect products
				with applications to continuum theories,
				{\it Adv. in Math.}, {\bf 137}, (1998), 1--81.

		\bibitem{MM} P.W. Michor, D. Mumford, \emph{An overview of the Riemannian metrics
				on spaces of curves using the Hamiltonian approach}, Appl. Comput.
				Harmon. Anal. {\bf 23} (1) (2007), 74--113.

\end{thebibliography}
\end{document}